\documentclass[11pt]{article}

\usepackage{amsmath,amssymb}

\topskip   \newtheorem{thm}{Theorem}
 \newtheorem{conj}[thm]{Conjecture}

 \newtheorem{prop}[thm]{Proposition}

 \newenvironment{proof}{\medskip\noindent{\it Proof.\ }}{\hfill \mbox{$\Box$}\medskip}
\title{A comment on Ryser's conjecture for intersecting hypergraphs}

\author{
Toufik Mansour \\
\small Department of Mathematics, University of Haifa, Haifa 31905, Israel\\[-.8ex]
\small{E-mail:toufik@math.haifa.c.il}\\[1.5ex]
Chunwei Song  \\
\small LMAM, School of Mathematical Sciences, Peking University, 1000871 Beijing, P.R. China \\[-0.8ex]
\small {E-mail: csong@math.pku.edu.cn}\\[1.5ex]
Raphael Yuster \\
\small Department of Mathematics, University of Haifa, Haifa 31905, Israel\\[-.8ex]
\small{E-mail:raphy@math.haifa.c.il}
}
\date{}

\begin{document}
\maketitle

\begin{abstract}
Let $\tau(\mathcal{H})$ be the cover number and
$\nu(\mathcal{H})$ be the matching number of a hypergraph
$\mathcal{H}$. Ryser conjectured that every
$r$-partite hypergraph $\mathcal{H}$ satisfies the inequality
$\tau(\mathcal{H}) \leq (r-1) \nu (\mathcal{H})$.
This conjecture is open for all $r \ge 4$.
For intersecting hypergraphs, namely those with $\nu(\mathcal{H})=1$,
Ryser's conjecture reduces to $\tau(\mathcal{H}) \leq r-1$.
Even this conjecture is extremely difficult and is open for all $ r \ge 6$.
For infinitely many $r$ there are examples of intersecting $r$-partite hypergraphs
with $\tau(\mathcal{H})=r-1$, demonstrating the tightness of the conjecture for such $r$.
However, all previously known constructions are not optimal as
they use far too many edges. How sparse can an intersecting $r$-partite hypergraph
be, given that its cover number is as large as possible, namely $\tau(\mathcal{H}) \ge r-1$?
In this paper we solve this question for $r \le 5$, give an almost optimal
construction for $r=6$, prove that any $r$-partite intersecting hypergraph
with $\tau(H) \ge r-1$ must have at least
$(3-\frac{1}{\sqrt{18}})r(1-o(1)) \approx 2.764r(1-o(1))$ edges,
and conjecture that there exist constructions with $\Theta(r)$ edges.
\bigskip

\noindent \emph{Keywords:} hypergraph, Ryser's conjecture, covering, matching, $r$-partite, intersecting.

\noindent \small  \emph{Mathematics Subject Classification}:
05C65, 05D05, 05C75
\end{abstract}

%%% ----------------------------------------------------------------------

\section{Introduction}
For a hypergraph $\mathcal{H}=(V, E)$, the (vertex) {\em cover number},
denoted by $\tau(\mathcal{H})$, is the minimum size of a
vertex set that intersects every edge. The {\em matching
number}, denoted by $\nu(\mathcal{H})$, is the maximum size
of a subset of edges whose elements are pairwise-disjoint.

Clearly, $\tau(\mathcal{H}) \le \nu(\mathcal{H})$ for any hypergraph.
In the graph-theoretic case, K\"onig's Theorem \cite{Ko-1936} asserts that the converse non-trivial inequality also holds for bipartite graphs. Thus, if ${\mathcal H}$ is a bipartite graph then $\tau(\mathcal{H}) = \nu(\mathcal{H})$.
Ryser conjectured the following hypergraph generalization of K\"onig's Theorem for hypergraphs.
A hypergraph is called {\em $r$-partite} if its vertex set can be partitioned into $r$ parts, and every edge
contains precisely one vertex from each part. In particular, $r$-partite hypergraphs are $r$-uniform.
Ryser conjectured that every $r$-partite hypergraph $\cal{H}$ satisfies
$\nu(\mathcal{H}) \le (r-1)\tau(\mathcal{H})$.
This conjecture turns out to be notoriously difficult. Indeed only the case $r=3$ has been proved
by Aharoni \cite{Ah-2001} using topological methods.

A hypergraph is called {\em intersecting} if any two edges have nonempty intersection.
Clearly, ${\mathcal H}$ is intersecting if and only if $\nu(\mathcal{H})=1$.
For intersecting hypergraphs, Ryser's conjecture amounts to:
\begin{conj}
\label{conj:Ryser2}
If ${\mathcal H}$ is an $r$-partite intersecting hypergraph then $\tau(\mathcal{H}) \leq r-1$.
\end{conj}
Conjecture \ref{conj:Ryser2} is still wide open.
It has been proved for $r=4,5$ by Tuza \cite{Tu-1979,Tu-1983}.
We note that the case $r=3$ of Conjecture \ref{conj:Ryser2} was first proved by Tuza in \cite{Tu-1983},
before Aharoni's general proof for the case $r=3$.

Conjecture \ref{conj:Ryser2} (if true) is tight in the sense that for infinitely many $r$
there are constructions of intersecting $r$-partite hypergraphs with $\tau(\mathcal{H}) = r-1$. Indeed, whenever $r=q+1$ and $q$ is a prime power,
consider the finite projective plane of order $q$ as a hypergraph.
This hypergraph is $r$-uniform and intersecting. To make it $r$-partite one just needs to delete one point from the projective plane. This truncated projective plane
gives an intersecting $r$-partite hypergraph with cover number $r-1$, with
$q^2+q=r(r-1)$ vertices, and with $q^2=(r-1)^2$ edges.

However, the projective plane construction is not the ``correct''
extremal construction, not only because it does not apply to all
$r$, but also because it is not the smallest possible. Although the
projective plane construction only contains $r(r-1)$ vertices (and
this is clearly optimal since otherwise some vertex class would have
size less than $r-1$ resulting in a cover number less than
$r-1$), should an extremal example contain so many (namely,
$(r-1)^2$) edges? In order to understand the extremal behavior of
intersecting $r$-partite hypergraphs, it is desirable to construct
the {\em sparsest} possible intersecting $r$-partite hypergraph with
cover number as {\em large} as possible, namely at least $r-1$.

More formally, let $f(r)$ be the minimum integer so that there exists an $r$-partite intersecting hypergraph $\mathcal{H}$ with $\tau(\mathcal{H}) \ge r-1$ and with $f(r)$ edges.
(we write $\tau(\mathcal{H}) \ge r-1$ instead of $\tau(\mathcal{H}) = r-1$ to allow for the possibility that Conjecture \ref{conj:Ryser2} is false; also note that trivially $\tau(\mathcal{H}) \le r$ since the set of vertices of any edge forms a cover).
A trivial lower bound for $f(r)$ is $2r-3$. Indeed, the edges of an intersecting hypergraph with at most $2r-4$ edges can greedily be covered with $r-2$ vertices.
We prove, however, the following non-trivial lower bound.
\begin{thm}
\label{t-lower}
$f(r) \ge (3-\frac{1}{\sqrt{18}})r(1-o(1)) \approx 2.764r(1-o(1))$.
\end{thm}

Although we do not have a matching upper bound, we conjecture that a linear (in $r$) number of edges indeed suffice.
\begin{conj}
\label{conj:sparse}
$f(r) = \Theta(r)$.
\end{conj}

Computing precise values of $f(r)$ seems to be a difficult problem.
Trivially, $f(2)=1$. It is also easy to see that $f(3)=3$.
Indeed, a $3$-partite intersecting hypergraph with only two edges has cover number $1$. The hypergraph whose edges are $(a_1,b_1,c_1), (a_1,b_2,c_2),(a_2, b_1, c_2)$
is a $3$-partite intersecting hypergraph with cover number $2$.
The next theorem establishes the first non-trivial values of $f(r)$,
namely $r=4,5$, in addition to upper and lower bounds in the case $r=6$.
More specifically, we prove:
\begin{thm}
\label{t1} $f(4)=6$, $f(5)=9$, and $12 \le f(6) \le 15$.
\end{thm}
Comparing our constructions with the projective plane construction, we see that
in the case $r=4,5,6$ the latter has $9,16$ and $25$ edges respectively.
Thus, the projective plane construction is far from being optimal.
Our constructions also have the property that the number of vertices they contain is
$r(r-1)$, which, as mentioned earlier, is optimal.

In the rest of this paper we prove Theorems \ref{t-lower} and \ref{t1}.

%%% ----------------------------------------------------------------------

\section{Proof of Theorem \ref{t-lower}}

Throughout this section we assume that ${\mathcal H}$ is an
$r$-partite intersecting hypergraph with $\tau(\mathcal{H}) \ge
r-1$. Recall that the {\em degree} of a vertex $v$ in
a hypergraph is the number of edges containing $v$.

Consider the following greedy procedure, starting with the original hypergraph
${\mathcal H}$. As long as there is a vertex $x$ of degree at least $4$
in the current hypergraph, we delete $x$ and all of the edges containing $x$
from the current hypergraph, thereby obtaining a smaller hypergraph.
Vertices that become isolated are also deleted.
Denote by ${\mathcal H}_3$ the hypergraph obtained at the end of the greedy procedure
and denote by $X_4$ the set of vertices deleted by the greedy procedure.
Notice that ${\mathcal H}_3$ is either the empty hypergraph or else it is an $r$-partite intersecting hypergraph, every vertex of which has degree at most $3$.
We then continue in the same manner, where as long as there is a vertex $x$ of degree $3$
in the current hypergraph, we delete $x$ and all of the edges containing $x$
from the current hypergraph. Again, vertices that become isolated are also deleted.
Denote by ${\mathcal H}_2$ the hypergraph obtained at the end of this second greedy procedure
and denote by $X_3$ the set of vertices deleted in the second greedy procedure.
Notice that ${\mathcal H}_2$ is either the empty hypergraph or else it is an $r$-partite intersecting hypergraph, every vertex of which has degree at most $2$.

We first claim that ${\mathcal H}_3$ contains at most $2r+1$ edges.
Indeed, if $H$ is any edge, then every vertex of $H$ appears in at most
two other edges. Thus, there are at most $2r$ other edges in addition to $H$.
Similarly, ${\mathcal H}_2$ contains at most $r+1$ edges.
Let, therefore, the number of ${\mathcal H}_3$ be denoted by $\gamma r$
and hence $0 \le \gamma \le 2+1/r$.

Consider first the case $\gamma \le 1$.
In this case we can cover the edges of ${\mathcal H}_3$ greedily with
a set $U$ of at most $\lceil \gamma r /2 \rceil$ vertices. Now,
since $U \cup X_4$ is a cover of ${\mathcal H}$ and since
$\tau({\mathcal H}) \ge r-1$, we have that
$|X_4| \ge r-1 - \lceil \gamma r /2 \rceil$.
As every vertex of $X_4$ was greedily selected to appear in four {\em distinct} edges
of ${\mathcal H} - {\mathcal H}_3$ we have that the number of edges of ${\mathcal H}$
is at least
$$
4|X_4|+\gamma r \ge 4(r-1 - \lceil \gamma r /2 \rceil)+ \gamma r = (4-\gamma)r-6
\ge 3r-6
$$
which is even better than the bound in the statement of the theorem.

We may now assume that $1 < \gamma \le 2+1/r$.
Since ${\mathcal H_2}$ has at most $r+1$ edges, we have that
$|X_3| \ge (\gamma r - r - 1)/3$. The number of edges of ${\mathcal H}_2$
is $\gamma r - 3|X_3|$. It follows that there is a cover of
${\mathcal H}_3$ whose size is at most
$$
|X_3| + \left\lceil \frac{\gamma r-3|X_3|}{2} \right\rceil \le \left(\frac{1}{6}+\frac{\gamma}{3}\right)r(1+o(1)).
$$
As such a cover, together with $X_4$, is a cover of ${\mathcal H}$,
and since $\tau(\mathcal{H}) \ge r-1$, we have that
$$
|X_4| \ge \left(\frac{5}{6}-\frac{\gamma}{3}\right)r(1-o(1)).
$$
We therefore have that the number of edges of ${\mathcal H}$, which is at least $4|X_4|+\gamma r$, is at least
\begin{equation}
\label{e1}
\frac{10-\gamma}{3}r(1-o(1)).
\end{equation}

There is, however, another way to bound from below the number of edges of
${\mathcal H}$.  For $i=1,2,3$ let $\alpha_ir^2$ denote the number of vertices of ${\mathcal H}_3$ having degree $i$.
Since the sum of the degrees in ${\mathcal H}_3$ is
$\gamma r^2$ we have:
$$
\alpha_1 + 2\alpha_2 + 3\alpha_3 = \gamma.
$$
Consider a specific edge $H$  of  ${\mathcal H}_3$ and let $r\beta_i^H$ be the number of vertices in $H$ with degree $i$ for $i=1,2,3$.
Clearly, $\beta_1^H+\beta_2^H+\beta_3^H=1$.
As $H$ intersects every edge we must have
$$
r\beta_2^H + 2r\beta_3^H \ge \gamma r -1.
$$
It follows that:
$$
2\beta_1^H + \beta_2^H = 2 - \beta_2^H - 2\beta_3^H \le 2 - \gamma + 1/r.
$$
In particular,
$$
\sum_{H \in {\mathcal H}_3} (2\beta_1^H + \beta_2^H) \le \gamma r (2-\gamma+1/r).
$$
On the other hand, by definition we have that
$$
\sum_{H \in {\mathcal H}_3} r\beta_1^H = \alpha_1 r^2~,~ \qquad
\sum_{H \in {\mathcal H}_3} r\beta_2^H = 2\alpha_2 r^2.
$$
It follows that
$$
2\alpha_1 + 2\alpha_2 \le \gamma(2-\gamma+1/r).
$$
Hence,
$$
\alpha_1 + \alpha_2 + \alpha_3 =
\frac{\gamma}{3}+\frac{2}{3}\alpha_1 + \frac{1}{3} \alpha_2 \le \frac{\gamma}{3}+\frac{2}{3}\alpha_1 + \frac{2}{3} \alpha_2 \le
$$
$$
\frac{\gamma}{3}+\frac{1}{3}\gamma(2-\gamma+1/r) =
\gamma - \frac{1}{3}\gamma^2 + \frac{\gamma}{3r}.
$$
Since $r^2(\alpha_1+\alpha_2+\alpha_3)$ is the number of vertices of
${\mathcal H}_3$ we have that ${\mathcal H}_3$ has at most
$$
(\gamma - \frac{1}{3}\gamma^2)r^2(1+o(1))
$$
vertices. In particular, there is a vertex class consisting of at most
$$
(\gamma - \frac{1}{3}\gamma^2)r(1+o(1))
$$
vertices. Since any vertex class of ${\mathcal H}_3$,
together with $X_4$, form a cover of ${\mathcal H}$, we have that
$$
|X_4| \ge (1-\gamma+\frac{1}{3}\gamma^2)r(1-o(1)).
$$
We therefore have that the number of edges of ${\mathcal H}$, which is at least $4|X_4|+\gamma r$, is at least
\begin{equation}
\label{e2}
(4-3\gamma+\frac{4}{3}\gamma^2)r(1-o(1)).
\end{equation}
Comparing (\ref{e1}) and (\ref{e2}) we see that the minimum of the maximum of both
of them is attained when $\gamma = 1 + 1/\sqrt{2}$, and in this case the number of
edges of ${\mathcal H}$ is at least $(3-\frac{1}{\sqrt{18}})r(1-o(1))$, as
required.
{\hfill \mbox{$\Box$}\medskip}

%%% ----------------------------------------------------------------------

\section{Proof of Theorem \ref{t1}}

\subsection{The case $r=4$}

We need to show first that $f(4) > 5$. Assume the contrary and let
${\mathcal H}$ be a $4$-partite intersecting hypergraph with only $5$ edges and
with $\tau(\mathcal{H}) \ge 3$. No vertex can appear in three or more edges,
since such a vertex $v$, and a vertex $u$ intersecting the (at most two) edges
in which $v$ does not appear form a cover of size $2$, a contradiction.
Thus, every vertex has degree at most $2$. Now, since there are $\binom{5}{2}$
nonempty intersections of pairs of edges of $\mathcal{H}$, we have, by the inclusion-exclusion principle that $\mathcal{H}$ contains at most
$5 \cdot 4 - \binom{5}{2} = 10$ vertices. But this means that some vertex class
contains at most two vertices, again resulting in $\tau(\mathcal{H}) \le 2$, a contradiction.

We construct a $4$-partite intersecting hypergraph with $6$ edges and
with $\tau(\mathcal{H}) = 3$.
Consider the four vertex classes $V_1=\{a_1,a_2,a_3\}$,
$V_2=\{b_1,b_2,b_3\}$, $V_3=\{c_1,c_2,c_3\}$, and $V_4=\{d_1,d_2,d_3\}$.
The $6$ edges are $(a_1, b_1, c_1, d_1)$, $(a_1, b_2, c_2, d_2)$,
$(a_2,b_1,c_2,d_3)$, $(a_2,b_2,c_3,d_1)$, $(a_3,b_3,c_2,d_1)$, and
$(a_3,b_1,c_3,d_2)$. It is easy to check that any two edges intersect
and that two vertices cannot cover all $6$ edges.

\subsection{The case $r=5$}

We need to show first that $f(5) > 8$. Assume the contrary and let
${\mathcal H}$ be a $5$-partite intersecting hypergraph with only $8$ edges and
with $\tau(\mathcal{H}) \ge 4$.

Notice first that there is no vertex with degree $4$ or greater, since if $v$ is such
a vertex, then the (at most) four remaining edges not containing $v$ can always
be greedily covered with two additional vertices, resulting in cover number at most $3$, a contradiction. Thus, we may assume that the degree of each vertex is
at most $3$. Clearly we can assume that ${\mathcal H}$ has at least $20$ vertices,
as otherwise there is a vertex class with at most three vertices, again contradicting the assumption that $\tau(\mathcal{H}) \ge 4$.

Let $x_i$ denote the number of vertices with degree $i$, for $i=1,2,3$.
Since the sum of the degrees of all vertices is $8 \cdot 5 = 40$ we have
that $x_1+2x_2+3x_3=40$.
We claim also that each vertex class has at most one vertex with degree $3$.
Indeed, if there were two such vertices in the same vertex class, then they both cover
$6$ edges, and the remaining two edges can be covered by an additional vertex,
contradicting the assumption that $\tau(\mathcal{H}) \ge 4$.
Thus, we have that $x_3 \le 5$.

As there are $28=\binom{8}{2}$ intersections, we have that $x_2+3x_3
\ge 28$. Finally, since there are at least $20$ vertices, we have
$x_1+x_2+x_3 \ge 20$. Now, it follows that $x_1+2x_2+4x_3 \ge 48$
which implies that $x_3 \ge 8$ which contradicts $x_3 \le 5$.

\vskip0.1cm

Next we construct a $5$-partite intersecting hypergraph with $9$
edges and with $\tau(\mathcal{H}) = 4$. Consider the five vertex
classes $V_1=\{1,2,3,4\}$, $V_2=\{5,6,7,8\}$, $V_3=\{9,10,11,12\}$,
$V_4=\{13,14,15,16\}$, and $V_5=\{17,18,19,20\}$. The $9$ edges are
divided into two parts:
$$
A = \{(1,5,9,13,17), (2,6,10,14,17), (3,7,10,13,18), (1,6,11,15,18), (2,7,9,15,19)\},
$$
$$
B = \{(4,5,10,15,20), (4,7,11,16,17), (4,8,9,14,18), (4,6,12,13,19)\}.
$$
First, notice that the constructed hypergraph is, indeed, $5$-partite, and intersecting. Also note that $\tau(\mathcal{H}) \le 4$ by considering, for example,
the cover $\{4,17,18,19\}$. It remains to show that that there is no cover of size $3$. There is only one vertex with degree $4$, and it is vertex $4$. Vertex $4$ precisely covers the set $B$. Notice, however, that
no vertex covers three edges of $A$. This means that any cover containing $4$ must have size at least $4$. We now only need to rule out the possibility of a
cover of size $3$, each vertex of which has degree $3$.
The vertices of degree $3$ are $\{1,2,6,7,9,10,13,15,17,18\}$. However, each of them
appears at most one time in $B$, hence any three of them cannot cover all the vertices of $B$.

\subsection{The case $r=6$}

We  first show that $f(6) > 11$.  Like before, assume the contrary
and let ${\mathcal H}$ be a $6$-partite intersecting hypergraph with
only $11$ edges and with $\tau(\mathcal{H}) \ge 5$.

Notice first that there is no vertex with degree $5$ or greater,
since if $v$ is such a vertex, then the (at most) $6$ remaining edges
not containing $v$ can always be greedily covered with $3$ additional
vertices, resulting in cover number at most $4$, a
contradiction. Thus, we may assume that the degree of each vertex is
at most $4$. Clearly we can assume that ${\mathcal H}$ has at least
$30$ vertices, as otherwise there is a vertex class with at most $4$
vertices, again contradicting the assumption that $\tau(\mathcal{H})
\ge 5$.

Let $x_i$ denote the number of vertices with degree $i$, for
$i=1,2,3,4$. Since the sum of the degrees of all vertices is $11 \cdot
6 = 66$ we have that $x_1+2x_2+3x_3+4x_4=66$.
Again, notice that the assumption that $\tau(\mathcal{H}) \ge 5$ implies that
each vertex class has at most one vertex with degree $4$,
at most two vertices of degree $3$,
and if  vertex class contains a vertex of degree $4$, it does not contain
a vertex of degree $3$.
Thus, $x_4 \le 6$ and $x_3 \le 12-2x_4$. Hence $x_3+3x_4 \le 18$.

As there are $55=\binom{11}{2}$ intersections, we have that $x_2+3
x_3 +6x_4 \ge 55$.  Combining this with the fact $x_1+x_2+x_3+x_4 \ge 30$ we
have  that $x_1+2x_2+4x_3 + 7x_4 \ge 85$.  This implies that $x_3+3x_4 \ge 19$,
contradicting the fact that $x_3+3x_4 \le 18$.

\vskip0.1cm

Next we create a $6$-partite intersecting hypergraph $\mathcal{H}=(V,E)$
with $30$ vertices and $15$ edges as follows.
The six vertex classes of $V$ are:
$$
\begin{array}{lll}
V_1=\{a_1,a_3,a_4,a_6,a_8\},&
V_2=\{b_1,b_2,b_4,b_8,b_{12}\},&
V_3=\{c_1,c_2,c_4,c_7,c_{11}\},\\
V_4=\{d_1,d_2,d_3,d_9,d_{11}\},&
V_5=\{e_1,e_2,e_3,e_5,e_7\},&
V_6=\{f_1,f_2,f_3,f_5,f_7\}.
\end{array}
$$
(The fact that the selection of indices in each set is not consecutive simplifies the
description that follows.) We construct $E$ in several steps so as to
guarantee that
\begin{enumerate}
  \item $| H \cap V_i | = 1 , \forall H \in E, 1 \leq i \leq 6$, so that $\mathcal{H}$ is a $6$-partite.
  \item $H \cap H' \neq \emptyset, \forall H, H' \in E$.
  \item $\tau(\mathcal{H}) \geq 5$.
\end{enumerate}

\textbf{Step 1}: A ``cyclic" construction.  Throughout the whole
procedure we try to let the edges, albeit intersecting, repeat as
little as possible.
$$\begin{array}{lll}
H_1 = \{ a_1, b_1, c_1, d_1, e_1, f_1 \}, & H_2 = \{ a_1, b_2, c_2, d_2, e_2, f_2 \}, & H_3 = \{ a_3, b_1, c_2, d_3, e_3, f_3 \},\\
H_4 = \{ a_4, b_4, c_4, d_1, e_2, f_3 \}, & H_5 = \{ a_3, b_4,
c_1, d_2, e_5, f_5 \}, & H_6 = \{ a_6, b_2, c_4, d_3, e_5, f_1 \}.
\end{array}$$
Note that by construction,
$$\begin{array}{ll}
| H_i \cap H_j | =1,& 1 \leq i < j \leq 6,\\
H_i \cap H_j \cap H_k = \emptyset,& 1 \leq i < j <k \leq 6.
\end{array}$$
Therefore, a minimum cover of $H_1, \ldots, H_6$ has size
$3$. Now we consider the pairwise intersections and take the
union of every three mutually disjoint pairs (the union of the
three pairs forms exactly $\bigcup_{i=1}^{6} H_i$), for instance
$\{ H_1 \cap H_2, H_3 \cap H_4, H_5 \cap H_6 \} = \{ a_1, f_3, e_5
\} $. The following list $\mathcal{L}_1$ thus contains all the
minimum covers of $H_1, \ldots, H_6$ :
\begin{align}
\mathcal{L}_1 = \bigl\{ & \{a_1, f_3, e_5\}, \{a_1, a_3, c_4\},
\{a_1,
d_3, b_4\},  \{b_1, e_2, e_5\}, \{b_1, d_2, c_4\}, \notag \\
& \{b_1, b_2, b_4\}, \{d_1, c_2, e_5 \}, \{d_1, d_2, d_3\}, \{d_1,
b_2, a_3\}, \{c_1, c_2, c_4\},  \notag
\\ & \{c_1, e_2, d_3 \}, \{c_1, b_2, f_3\},  \{ f_1, c_2, b_4 \}, \{ f_1, e_2, a_3 \}, \{ f_1, d_2, f_3 \} \bigr\}. \notag
\end{align}

\textbf{Step 2}: Any additional edge must contain a cover of $H_1,
\ldots, H_6$.  Selecting carefully an element from $\mathcal{L}_1$
each time, we construct the edges $H_7$ through $H_{10}$.
$$\begin{array}{ll}
H_7 = \{ a_1, b_4, c_7, d_3, e_7, f_7 \},&H_8 = \{ a_8, b_8, c_2, d_1, e_5, f_7 \}, \\
H_9 = \{ a_8, b_2, c_1, d_9, e_7, f_3 \}, & H_{10} = \{ a_3, b_8,
c_7, d_9, e_2, f_1 \}.
\end{array}$$
Up to now, $| H_i \cap H_j | =1, for 1 \leq i < j \leq 10$.  Moreover,
since the (unique) intersecting element of any pair of edges
constructed in Step 2 has a subscript index in $\{ 7, 8, 9\}$,
$H_i \cap H_j \cap H_k = \emptyset, 1\leq i<j<k \leq 10, \text{if either} \ |\{i, j, k \} \cap [6] | =3 \
 \text{or} \ |\{i, j, k \} \cap \{ 7, 8, 9, 10\} | \geq 2$.
Notice also that any minimum cover that covers $H_1-H_{10}$ consists of at least
four vertices.

\noindent
\textbf{Step 3}: Five additional edges that force an increase in
the size of the minimum cover.
$$
\begin{array}{lll} H_{11} = \{ a_1, b_8, c_{11},
d_{11}, e_5, f_3 \}, & H_{12} = \{ a_8, b_{12}, c_1, d_3, e_2, f_3
\}, & H_{13} = \{ a_8, b_4, c_2,
d_{11}, e_{1}, f_1 \}, \\
H_{14} = \{ a_{4}, b_8, c_{1}, d_3, e_2, f_1 \},& H_{15} = \{ a_8,
b_{8}, c_1, d_3, e_2, f_5 \}.
\end{array}
$$
Notice that the $15$ constructed edges indeed form an intersecting $6$-partite
hypergraph. It remains to show that:
\begin{prop}
\label{thm:r=6}
$\tau(\mathcal{H}) = 5$.
\end{prop}
\begin{proof}
Suppose the proposition is false, and
let $C =\{x, y, z, w\}$ be a cover of size $4$.  We use a
sequence of arguments to deduce that this is impossible.

For convenience, for $7 \leq i \leq 15$, let ${H_i}^*$ be the
collection of those vertices of $H_i$ with index subscript no
bigger than $6$. That is,
$$
\begin{array}{lll}
H_{7}^* = \{ a_1, b_4, d_3 \}, & H_{8}^* = \{ c_2, d_1, e_5 \}, &
H_{9}^* = \{ b_2, c_1, f_3 \},\\
H_{10}^* = \{ a_3, e_2, f_1 \},& H_{11}^* = \{ a_1, e_5, f_3 \}, &
H_{12}^* = \{ c_1, d_3, e_2, f_3 \},\\
H_{13}^* = \{ b_4, c_2, e_1, f_1\}, & H_{14}^* = \{ a_4, c_1, d_3,
e_2, f_1 \}, & H_{15}^* = \{ c_1, d_3, e_2, f_5 \}.
\end{array}
$$

\textbf{First}, notice that if $C \cap \{ b_{12}, c_{11}, d_{11} \}
\neq \emptyset$, then $| C \backslash \{ b_{12}, c_{11}, d_{11} \} |
\leq 3$.  So $ C \backslash \{ b_{12}, c_{11}, d_{11} \}$ must be a
triple, and furthermore it must be one of the triples in
$\mathcal{L}_1$. But no element of $\mathcal{L}_1$ may cover $H_7,
H_8, H_9$ and $H_{10}$ since $H_{7}^* - H_{10}^*$ are pairwise
disjoint.  Hence, $C \cap \{ b_{12}, c_{11}, d_{11} \} = \emptyset$.

\textbf{Second}, assume $C \cap \{a_{8}, b_{8} \} \neq \emptyset$.
Without loss of generality, let $x = C \cap \{a_{8}, b_{8} \}$.  Then $\{ y, z, w \}
\in \mathcal{L}_1$.

\begin{itemize}
\item Case 1: $x=a_8$. The fact that $\{ y, z, w \}$ must cover $H_7^*$
implies that $\{ y, z, w \} \cap \{ a_1, b_4, d_3\} \neq \emptyset$;
$\{ y, z, w \}$ must cover $H_{10}^*$ implies that $\{ y, z, w \}
\cap \{ a_3, e_2, f_1\} \neq \emptyset$; $\{ y, z, w \}$ must cover
$H_{11}^*$ implies that $\{ y, z, w \} \cap \{ a_1, e_5, f_3\} \neq
\emptyset$.  The only triple in $\mathcal{L}_1$ satisfying these
three requirements is $\{ a_1, a_3, c_4\}$.  But then $H_{14}^* =\{
a_4, c_1, d_3, e_2, f_1 \}$ is left uncovered. Impossible.
  \item Case 2: $x=b_8$.  The fact that $\{ y, z, w \}$ must cover $H_7^*$
implies that $\{ y, z, w \} \cap \{ a_1, b_4, d_3\} \neq \emptyset$;
$\{ y, z, w \}$ must cover $H_{9}^*$ implies that $\{ y, z, w \}
\cap \{ b_2, c_1, f_3\} \neq \emptyset$; $\{ y, z, w \}$ must cover
$H_{13}^*$ implies that $\{ y, z, w \} \cap \{ b_4, c_2, e_1, f_1\}
\neq \emptyset$.  The only triple in $\mathcal{L}_1$ satisfying
these three requirements is $\{ b_1, b_2, b_4\}$.  But then
$H_{12}^* =\{ c_1, d_3, e_2, f_3 \}$ is left uncovered. Impossible.
\end{itemize}
Hence, $C \cap \{a_{8}, b_{8} \} = \emptyset$.

\textbf{Third}, assume $C \cap \{c_{7}, e_7, f_7, d_9 \} \neq
\emptyset$. Let $C \cap \{c_{7}, e_7, f_7, d_9 \} =x$.
Then $\{ y, z, w \} $ must be a triple in $\mathcal{L}_1$, and it
must also cover $H_{11}^* - H_{15}^*$.  The fact that $\{ y, z, w
\}$ must cover $H_{11}^*$ implies that $\{ y, z, w \} \cap \{ a_1,
e_5, f_3\} \neq \emptyset$; $\{ y, z, w \}$ must cover $H_{13}^*$
implies that $\{ y, z, w \} \cap \{ b_4, c_2, e_1, f_1\} \neq
\emptyset$; $\{ y, z, w \}$ must cover $H_{15}^*$ implies that $\{
y, z, w \} \cap \{ c_1, d_3, e_2, f_5\} \neq \emptyset$.  The only
triple in $\mathcal{L}_1$ satisfying these three requirements is $\{
a_1, b_4, d_3\}$.  But then $H_8, H_9, H_{10}$ are not yet covered.
This can not be fixed by any additional one vertex as
$H_8 \cap H_9 \cap H_{10} = \emptyset$. Hence we must have
$C \cap \{c_{7}, e_7, f_7, d_9 \} = \emptyset$.

\textbf{Last}, we have by now established that the index of any
vertex in $C$ is in $\{ 1, 2,3,4,5,6 \}$.  Furthermore, $|C
\cap H_{i}^*|=1$, for $7 \leq i \leq 10$.  In particular, let $x= C \cap
H_{7}^*$, we discuss each of the three possibilities.

(i) $x=a_1$. Then $\{ y, z, w \}$ must cover $H_3 - H_6, H_8^* -
H_{10}^*, H_{12}^* - H_{15}^*$. Let $y \in H_8^*$.
\begin{itemize}
\item Case 1: $y=c_2$.  The fact $\{z, w \}$ needs to cover $H_4 -
H_6$ $\Rightarrow$ $\{z, w \} \cap \{ b_4, e_5, c_4\} \neq
\emptyset$;
 $\{z, w \}$ needs to cover $H_9^*$ $\Rightarrow$ $\{z, w \}
\cap \{ b_2, c_1, f_3\} \neq \emptyset$;  $\{z, w \}$ needs to cover
$H_{10}^*$ $\Rightarrow$ $\{z, w \} \cap \{ a_3, e_2, f_1\} \neq
\emptyset$.  But this is clearly impossible.

\item Case 2: $y=d_1$. The fact $\{z, w \}$ needs to cover $H_3,
H_5$ and $H_6$ $\Rightarrow$ $\{z, w \} \cap \{ a_3, e_5, d_3\}
\neq \emptyset$;
 $\{z, w \}$ needs to cover $H_9^*$ $\Rightarrow$ $\{z, w \}
\cap \{ b_2, c_1, f_3\} \neq \emptyset$;  $\{z, w \}$ needs to cover
$H_{13}^*$ $\Rightarrow$ $\{z, w \} \cap \{ b_4, c_2, e_1, f_1\}
\neq \emptyset$.  Impossible.

\item Case 3: $y=e_5$. The fact $\{z, w \}$ needs to cover $H_9^*$
$\Rightarrow$ $\{z, w \} \cap \{ b_2, c_1, f_3\} \neq \emptyset$;
$\{z, w \}$ needs to cover $H_{10}^*$ $\Rightarrow$ $\{z, w \} \cap
\{ a_3, e_2, f_1\} \neq \emptyset$; $\{z, w \}$ needs to cover
$H_{13}^*$ $\Rightarrow$ $\{z, w \} \cap \{ b_4, c_2, e_1, f_1\}
\neq \emptyset$; $\{z, w \}$ needs to cover $H_{15}^*$ $\Rightarrow$
$\{z, w \} \cap \{ c_1, d_3, e_2, f_5 \} \neq \emptyset$.  Putting
the four requirements all together, $\{z, w \}$ is forced to equal
$\{ c_1, f_1\}$.  But then $H_4$ is left uncovered.  Impossible.
\end{itemize}

(ii) $x=b_4$. Then $\{ y, z, w \}$ must cover $H_1 - H_3, H_6,
H_8^* - H_{12}^*, H_{14}^*, H_{15}^*$. Let $y \in H_8^*$.
\begin{itemize}

\item Case 1: $y=c_2$.  The fact $\{z, w \}$ needs to cover
$H_9^*$ $\Rightarrow$ $\{z, w \} \cap \{ b_2, c_1, f_3\} \neq
\emptyset$;  $\{z, w \}$ needs to cover $H_{10}^*$ $\Rightarrow$
$\{z, w \} \cap \{ a_3, e_2, f_1\} \neq \emptyset$; $\{z, w \}$
needs to cover $H_{11}^*$ $\Rightarrow$ $\{z, w \} \cap \{ a_1,
e_5, f_3\} \neq \emptyset$; $\{z, w \}$ needs to cover $H_{15}^*$
$\Rightarrow$ $\{z, w \} \cap \{ c_1, d_3, e_2, f_5 \} \neq
\emptyset$.  The four requirements force $\{z, w \}$  to equal $\{
f_3, e_2\}$.  But then $H_1$ is left uncovered. Impossible.

\item Case 2: $y=d_1$. The fact $\{z, w \}$ needs to cover $H_2,
H_3$ and $H_6$ $\Rightarrow$ $\{z, w \} \cap \{ c_2, b_2, d_3\}
\neq \emptyset$; $\{z, w \}$ needs to cover $H_{10}^*$
$\Rightarrow$ $\{z, w \} \cap \{ a_3, e_2, f_1\} \neq \emptyset$;
$\{z, w \}$ needs to cover $H_{11}^*$ $\Rightarrow$ $\{z, w \}
\cap \{ a_1, e_5, f_3\} \neq \emptyset$.  Impossible.

  \item Case 3: $y=e_5$. The fact $\{z, w \}$ needs to cover $H_1 - H_3$
$\Rightarrow$ $\{z, w \} \cap \{ a_1, b_1, c_2\} \neq \emptyset$;
 $\{z, w \}$ needs to cover $H_9^*$ $\Rightarrow$ $\{z, w \}
\cap \{ b_2, c_1, f_3\} \neq \emptyset$;  $\{z, w \}$ needs to
cover $H_{10}^*$ $\Rightarrow$ $\{z, w \} \cap \{ a_3, e_2, f_1\}
\neq \emptyset$.  Impossible.
\end{itemize}

(iii) $x=d_3$. Then $\{ y, z, w \}$ must cover $H_1, H_2, H_4, H_5,
H_8^* - H_{11}^*, H_{13}^*$. Let $y \in H_8^*$.
\begin{itemize}
\item Case 1: $y=c_2$.  The fact $\{z, w \}$ needs to cover $H_1,
H_4$ and  $H_5$ $\Rightarrow$ $\{z, w \} \cap \{ d_1, b_4, c_1\}
\neq \emptyset$;
 $\{z, w \}$ needs to cover $H_9^*$ $\Rightarrow$ $\{z, w \}
\cap \{ b_2, c_1, f_3\} \neq \emptyset$;  $\{z, w \}$ needs to cover
$H_{10}^*$ $\Rightarrow$ $\{z, w \} \cap \{ a_3, e_2, f_1\} \neq
\emptyset$.  Impossible.

\item Case 2: $y=d_1$. The fact $\{z, w \}$ needs to cover $H_9^*$
$\Rightarrow$ $\{z, w \} \cap \{ b_2, c_1, f_3\} \neq \emptyset$;
$\{z, w \}$ needs to cover $H_{10}^*$ $\Rightarrow$ $\{z, w \} \cap
\{ a_3, e_2, f_1\} \neq \emptyset$; $\{z, w \}$ needs to cover
$H_{11}^*$ $\Rightarrow$ $\{z, w \} \cap \{ a_1, e_5, f_3\} \neq
\emptyset$; $\{z, w \}$ needs to cover $H_{13}^*$ $\Rightarrow$
$\{z, w \} \cap \{ b_4, c_2, e_1, f_1 \} \neq \emptyset$. The four
requirements force $\{z, w \}$ to be $\{ f_3, e_2\}$. But then $H_2$
is left uncovered. Impossible.

 \item Case 3: $y=e_5$. The fact $\{z, w \}$ needs to cover $H_1$, $H_2$ and $H_4$
$\Rightarrow$ $\{z, w \} \cap \{ a_1, d_1, e_2\} \neq \emptyset$;
 $\{z, w \}$ needs to cover $H_9^*$ $\Rightarrow$ $\{z, w \}
\cap \{ b_2, c_1, f_3\} \neq \emptyset$; $\{z, w \}$ needs to cover
$H_{13}^*$ $\Rightarrow$ $\{z, w \} \cap \{ b_4, c_2, e_1, f_1 \}
\neq \emptyset$.  Impossible.
\end{itemize}

In conclusion, our assumption is contradicted. Hence,
$\tau(\mathcal{H}) = 5$.
\end{proof}

%%% ----------------------------------------------------------------------

\end{document}